\documentclass[11pt]{article}

\usepackage{latexsym}
\usepackage{amssymb}

\newtheorem{theorem}{Theorem}
\newtheorem{proposition}{Proposition}

\newtheorem{lemma}{Lemma}

\begin{document}
{
\begin{center}
{\Large\bf
Positive orthogonalizing weights on the unit circle for the generalized Bessel polynomials.}
\end{center}
\begin{center}
{\bf Sergey M. Zagorodnyuk}
\end{center}

\noindent
\textbf{Abstract.}
In this paper we study the generalized Bessel polynomials $y_n(x,a,b)$ (in the notation of Krall and Frink).
Let $a>1$, $b\in(-1/3,1/3)\backslash\{ 0\}$.  
In this case we present the following positive continuous weights $p(\theta) = p(\theta,a,b)$ on the unit circle for $y_n(x,a,b)$:
$$ 2\pi p(\theta,a,b) = -1 + 2(a-1) \int_0^1 e^{-bu\cos\theta} \cos(bu\sin\theta) (1-u)^{a-2} du, $$
where $\theta\in[0,2\pi]$. Namely, we have
$$ \int_0^{2\pi} y_n(e^{i\theta},a,b) y_m(e^{i\theta},a,b) p(\theta,a,b) d\theta = C_n \delta_{n,m},\qquad C_n\not=0,\ n,m\in\mathbb{Z}_+. $$
Notice that this orthogonality differs from the usual orthogonality of OPUC.
Some applications of the above orthogonality are given.

\noindent
\textbf{Key words}: generalized Bessel polynomials, orthogonality, complex Jacobi matrix, orthogonal polynomials.

\noindent
\textbf{MSC 2020}: 42C05.

\section{Introduction.}

The theory of orthogonal polynomials, especially the cases of the real line (OPRL) and of the unit circle (OPUC), 
attracts a lot of attention of researchers and it has 
a lot of contributions and applications, see~\cite{cit_50000_Gabor_Szego}, \cite{cit_3000_Chihara}, \cite{cit_5000_Ismail}, 
\cite{cit_48000_Simon_1}, \cite{cit_48000_Simon_2}, \cite{cit_3100_Grosswald}, \cite{cit_48500_Srivastava}.
Special attention has been devoted to four classical systems of orthogonal polynomials: Jacobi, Laguerre, Hermite
and Bessel polynomials. They are eigenfunctions of second-order linear differential operators with polynomial coefficients
(not depending on the index $n$).
The basic properties of the (generalized) Bessel polynomials $y_n(x,a,b)$ (here we shall follow the notation of Krall and Frink 
from~\cite{cit_5100_Krall_Frink}) are less transparent then those of OPRL and OPUC. 
That provides an 
intriguing challenge for mathematicians.  Question: are $y_n(x,a,b)$ orthogonal on the unit circle $\mathbb{T}$
or on the real line $\mathbb{R}$? The answer is that they are orthogonal on the both sets but not in a usual manner.
This paper is mainly devoted to the orthogonality of the generalized Bessel polynomials.

A question of the orthogonality of polynomials
\begin{equation}
\label{f1_3}
y_n(x) = y_n(x,a,b) := {}_2 F_0 \left(-n,n+a-1;-;-\frac{x}{b} \right),
\end{equation}
started  in 1949 with the basic paper of Krall and Frink~\cite{cit_5100_Krall_Frink}. In their definition of $y_n(x,a,b)$ Krall and Frink assumed
that $a$ is not a negative integer or zero, and $b$ is not zero. In fact, Krall and Frink throughout their paper 
did not write explicitly if parameters $a$ and $b$ could take complex values or they were reals.
Among other basic results for polynomials $y_n(x,a,b)$, Krall and Frink proved that $y_n(x,a,b)$ are orthogonal on the unit circle:
\begin{equation}
\label{f1_5}
\int_\mathbb{T} y_n(z,a,b) y_m(z,a,b) \rho(z,a,b) dz = c_n \delta_{n,m},\qquad c_n\not=0,\ n,m\in\mathbb{Z}_+,
\end{equation}
where
\begin{equation}
\label{f1_6}
\rho(z) = \rho(z,a,b) = \frac{1}{2\pi i} \sum_{j=0}^\infty \frac{\Gamma(a)}{\Gamma(a+j-1)} \left(
-\frac{b}{z}
\right)^j,\qquad z\in\mathbb{C}\backslash\{ 0 \};
\end{equation}
and
\begin{equation}
\label{f1_6_3}
c_k := \frac{ (-1)^{k+1} b k! }{ (k+a-1) [k+a-2]_{k-1} },\quad k\in\mathbb{N};\quad c_0=-b.
\end{equation}
Here $[u]_0 := 1$, $[u]_k := u(u-1)...(u-k+1)$, $k\in\mathbb{N}$; $u\in\mathbb{C}$.
It is clear that in the case of ordinary Bessel polynomials ($a=b=2$) the weight function $\rho$ simplifies to the
following complex-valued weight:
$$ \rho(z,2,2)  = \frac{1}{2\pi i} e^{-\frac{2}{z}}. $$

In 1978 H.L.~Krall sent to Grosswald a handwritten manuscript describing real-valued orthogonalizing weights on the real line
for the generalized Bessel polynomials $y_n(x,a,b)$, see~\cite[p. 33]{cit_3100_Grosswald}.
In Grosswald's book~\cite{cit_3100_Grosswald} only statements are presented (without proofs).

Eight years later, in 1986, Exton presented the following orthogonality relations for Bessel polynomials~\cite{cit_3050_Exton} (see
also~\cite{cit_48500_Srivastava} for equivalent formulations):
$$ \frac{1}{2\pi i} \int_0^{(\infty+)} z^{a-2} e^{-1/z} y_n(z,a,1) y_m(z,a,1) dz = $$
\begin{equation}
\label{f1_7}
= (-1)^{a+n} \frac{n!}{ (a+2n-1)\Gamma(a+n-1) } \delta_{m,n},\qquad m,n\in\mathbb{Z}_+.
\end{equation}
Here the countour of integration is taken along a simple loop starting at the origin, encircling $\infty$ once in the positive direction,
and then returning to the origin.

In 1992, Kwon, Kim and Han presented two real-valued measures supported on $[0,+\infty)$ for
ordinary Bessel polynomials~\cite{cit_7100_Kwon}:
\begin{equation}
\label{f1_8}
w_1(x) = \left\{
\begin{array}{cc} 0, & x\leq 0,\\
 -e^{-2/x} \int_x^\infty t^{-2} e^{-t^{1/4} + 2/t} \sin t^{1/4} dt, & x>0;
\end{array}
\right.
\end{equation}
and
\begin{equation}
\label{f1_10}
w_2(x) = \left\{
\begin{array}{cc} 0, & x\leq 0,\\
 -e^{-2/x} \int_x^\infty t^{-2} e^{2/t} h(t) dt, & x>0.
\end{array}
\right.
\end{equation}
Here $h$ is due to Stieltjes:
$$ h(x) = \left\{
\begin{array}{cc} 0, & x\leq 0,\\
 \sin(2\pi\ln x) e^{-\ln^2 x}, & x>0.
\end{array}
\right.
$$

Evans, Everitt, Kwon and Littlejohn constructed real-valued orthogonalizing weights for the generalized Bessel polynomials for some values of
the parameters in~\cite{cit_3050_Evans__JCAM}.

About the same time, Duran presented two real-valued weights on $[0,+\infty)$ for the generalized Bessel polynomials~\cite{cit_3040_Duran}.
These weights involved some integrals with Bessel functions under the integral signs.

Finally, there were studied some quasi-orthogonality properties of Bessel polynomials, see~\cite{cit_48500_Srivastava} for details and
references.

Suppose that $a>1$, and $b\in\mathbb{R}\backslash\{ 0 \}$.
Consider a function $p(\theta)=p(\theta,a,b)$ defined as follows:
\begin{equation}
\label{f1_12}
2\pi p(\theta,a,b) = -1 + 2(a-1) \int_0^1 e^{-bu\cos\theta} \cos(bu\sin\theta) (1-u)^{a-2} du,\qquad \theta\in[0,2\pi].
\end{equation}
This function is positive if $|b|<1/3$, and a more precise numerical estimate will be given below.
In this case we have the orthogonality relations:
\begin{equation}
\label{f1_14}
\int_0^{2\pi} y_n(e^{i\theta},a,b) y_m(e^{i\theta},a,b) p(\theta,a,b) d\theta = - \frac{c_n}{b} \delta_{n,m},\quad  n,m\in\mathbb{Z}_+, 
\end{equation}
where $c_n$ are defined in~(\ref{f1_6_3}).

The function $p(\theta,a,b)$ can be expressed in terms of $\rho(z,a,b)$ in the following way:
$$ 2\pi p(\theta,a,b) = -1 + \frac{ 2(a-1) }{b} \cos\theta - $$
\begin{equation}
\label{f1_15}
- \frac{2\pi i}{b} \left(
e^{i\theta} \rho(e^{i\theta},a,b) + e^{-i\theta} \rho(e^{-i\theta},a,b)
\right),\qquad \theta\in[0,2\pi].
\end{equation}
This follows from the representation:
\begin{equation}
\label{f1_15_1}
2\pi i \rho(z,a,b) = a-1 - \frac{b}{z} {}_1 F_1 (1;a;-b/z),\qquad z\in\mathbb{C}\backslash\{ 0 \},
\end{equation}
and relations~(\ref{f2_22}),(\ref{f2_25}) below.

Denote by $C$ the operator of conjugation: $f(\theta) \mapsto \overline{f(\theta)}$ in the usual space $L^2_p$. Introduce
the $C$-form:
$$ [f,g]_C := (f,Cg)_{L^2_p},\qquad f,g\in L^2_p. $$
Then relation~(\ref{f1_14}) can be rewritten in the following form:
\begin{equation}
\label{f1_16}
[ y_n(e^{i\theta},a,b), y_m(e^{i\theta},a,b) ]_C = - \frac{c_n}{b} \delta_{n,m},\qquad n,m\in\mathbb{Z}_+. 
\end{equation}
This means that polynomials $y_n$ are $C$-orthogonal in the Hilbert space $L^2_p$.
Such kind of orthogonal systems were studied by Garcia and Putinar, see~\cite{cit_3070_Garcia}, \cite{cit_3080_Garcia_Putinar__Tohoku_Math_J}
and~\cite{cit_3090_Garcia_Putinar}. One of important questions here is the convergence of expansions of elements of the Hilbert space
into series by these systems.

On the other hand, Milovanovi\'c and his collegues studied orthogonal polynomials on a semicircle and on arcs of the unit circle
of type~(\ref{f1_14}), see a survey in~\cite{cit_7200_Milovanovic}.

We shall discuss some possible applications of the orthogonality relations~(\ref{f1_14}). We can give some
estimates for arbitrary solutions of the corresponding second-order difference equation. It is shown that 
for sufficiently small values of $b$
the use of the weight $p$ provides better estimates than the use of $\rho$.

\noindent
{\bf Notations. }
As usual, we denote by $\mathbb{R}, \mathbb{C}, \mathbb{N}, \mathbb{Z}, \mathbb{Z}_+$
the sets of real numbers, complex numbers, positive integers, integers and non-negative integers,
respectively; $\mathbb{T} := \{ z\in\mathbb{C}:\ |z|=1 \}$, $\mathbb{D} := \{ z\in\mathbb{C}:\ |z|<1 \}$. 
For $u\in\mathbb{C}$ we denote $(u)_0 := 1$, $(u)_k := u(u+1)...(u+k-1)$, $k\in\mathbb{N}$ (the shifted factorial or
Pochhammer's symbol).
For $u\in\mathbb{C}$ we also denote $[u]_0 := 1$, $[u]_k := u(u-1)...(u-k+1)$, $k\in\mathbb{N}$.
By $\mathbb{P}$ we mean a set of all complex polynomials.
By $\mathfrak{B}(M)$ we denote the set of all Borel subsets of a set $M\subseteq\mathbb{C}$.
For a measure $\mu$ on a $\sigma$-algebra $\mathfrak{A}$ of subsets of $\mathbb{C}$ we denote by $L^2_\mu$ the usual space of all
(classes of equivalence of) $(\mathfrak{A},\mathfrak{B}(\mathbb{C}))$ measurable complex-valued functions $f$ on $M$, such that
$\int_M |f|^2 d\mu < +\infty$.
If $\mu(\delta) = \int_\delta p(\theta) d\theta$, then we shall also use the notation $L^2_p = L^2_\mu$.

By $l^2$ we denote the usual space of square-summable complex sequences $\vec u = (u_k)_{k=0}^\infty$, $u_k\in\mathbb{C}$,
and $l^2_{fin}$ means the subset of all finitely supported vectors from $l^2$.

If $H$ is a Hilbert space then $(\cdot,\cdot)_H$ and $\| \cdot \|_H$ mean
the scalar product and the norm in $H$, respectively. 
Indices may be omitted in obvious cases.
For a linear operator $A$ in $H$, we denote by $D(A)$
its  domain, by $R(A)$ its range, and $A^*$ means the adjoint operator
if it exists. If $A$ is invertible then $A^{-1}$ means its
inverse. $\overline{A}$ means the closure of the operator, if the
operator is closable. If $A$ is bounded then $\| A \|$ denotes its
norm.

The generalized hypergeometric function is denoted by
$$ {}_m F_n(a_1,...,a_m; b_1,...,b_n;x) = \sum_{k=0}^\infty \frac{(a_1)_k ... (a_m)_k}{(b_1)_k ... (b_n)_k} \frac{x^k}{k!}, $$
where $m,n\in\mathbb{N}$, $a_j\in\mathbb{C}$, $b_l\in\mathbb{C}\backslash\{ 0,-1,-2,...\}$.
By $\Gamma(z)$ and $\mathrm{B}(z)$ we denote the gamma function and the beta function, respectively.

\section{Positive weights for Bessel polynomials.}

Let $a>1$, and $b\in\mathbb{R}\backslash\{ 0 \}$ be some fixed numbers.
Consider a function $p(\theta)=p(\theta,a,b)$ defined as in~(\ref{f1_12}).
Since
$$ \int_0^1 \left| e^{-bu\cos\theta} \cos(bu\sin\theta) (1-u)^{a-2} \right| du \leq
M \int_0^1 (1-u)^{a-2} du < \infty,\ M>0, $$
then $p(\theta,a,b)$ is well-defined.
It is clear that $p(\theta,a,b)$ is real-valued. We are especially interested in those values of the parameter $b$
for which $p(\theta,a,b)$ is positive. Denote
\begin{equation}
\label{f2_3}
I_2(\theta) = I_2(\theta,a,b) :=  
2(a-1) \int_0^1 e^{-bu\cos\theta} \cos(bu\sin\theta) (1-u)^{a-2} du,\ \theta\in[0,2\pi]. 
\end{equation}
Suppose that
\begin{equation}
\label{f2_5}
0 < |b| < R,\quad \mbox{for some $R$: } 0 < R \leq\pi/2. 
\end{equation}
Then $| b u \sin\theta | <\pi/2$, and therefore the function under the integral sign in~(\ref{f2_3}) is non-negative.
Since for $u\in[0,1], \theta\in[0,2\pi]$ we have
$$ e^{-bu\cos\theta} \geq e^{-|b|} \geq e^{-R},\quad \cos(bu\sin\theta)\geq \cos R, $$
then
$$ I_2 \geq 2(a-1) e^{-R} \cos R \int_0^1 (1-u)^{a-2} du = 2 e^{-R} \cos R\rightarrow 2, $$
as $R\rightarrow +0$.
This ensures that for sufficiently small values of $b$ the function $p(\theta,a,b)$ is positive.
Consider the following function:
$$ y_0(x) := 2 e^{-x} \cos x,\qquad x\in[0,\pi/2]. $$
Since $y_0'(x) = -2 e^{-x} (\cos x + \sin x) < 0$, then $y_0(x)$ is strictly decreasing on $[0,\pi/2]$.
Notice that $y_0(0) = 2$, $y_0(\pi/2) = 0$. Thus, there exists a unique value $x_0\in[0,\pi/2]$ such that
$y_0(x_0) = 2 e^{-x_0} \cos x_0 = 1$.
The equation
$$ 2\cos x = e^x $$
can be studied by numerical methods (e.g. by the bisection method). The function FindRoot of \textit{Mathematica} gives
the following approximation of $x_0$:
\begin{equation}
\label{f2_8}
x_0 \approx 0.539785.
\end{equation}
Consequently, if
\begin{equation}
\label{f2_10}
|b| < x_0, 
\end{equation}
then 
\begin{equation}
\label{f2_14}
p(\theta,a,b)>0,\qquad \theta\in[0,2\pi].
\end{equation}
We can also provide a simple upper bound for $|b|$ which ensures the validity of~(\ref{f2_14}).
We look for $\zeta\in[0,\pi/2]$ such that $y_0(\zeta)\geq 1$, which is equivalent to
\begin{equation}
\label{f2_16}
2\cos \zeta\geq e^\zeta \Leftrightarrow 4\cos^2 \zeta\geq e^{2\zeta}. 
\end{equation}
Observe that
$$ e^{2\zeta} = \sum_{k=0}^\infty \frac{1}{k!} (2\zeta)^k \leq \sum_{k=0}^\infty (2\zeta)^k = \frac{1}{1-2\zeta}, $$
if $\zeta< 1/2$. 
Since $\sin \zeta \leq \zeta$, then
$4 \cos^2\zeta \geq 4 - 4\zeta^2$.
Consequently, if $0<\zeta<1/2$ is such that
\begin{equation}
\label{f2_18}
4-4\zeta^2 \geq \frac{1}{1-2\zeta}, 
\end{equation}
then relation~(\ref{f2_16}) holds.
Since $4-4\zeta^2\geq 3$, then~(\ref{f2_18}) holds if
\begin{equation}
\label{f2_20}
\frac{1}{1-2\zeta}\leq 3  \Leftrightarrow \zeta\leq 1/3.
\end{equation}
So, if $|b|<1/3$ then relation~(\ref{f2_14}) holds.

\begin{theorem}
\label{t2_1}
Let $x_0$ be the unique number from the interval $[0,\pi/2]$ such that $2 \cos x_0 = e^{x_0}$.
Let $a>1$, $b\in(-x_0,x_0)\backslash\{ 0 \}\supseteq (-1/3,1/3)\backslash\{ 0 \}$, be fixed numbers.
The generalized Bessel polynomials $y_n(x,a,b)$ satisfy orthogonality relations~(\ref{f1_14}), where
the positive continuous weight $p(\theta) = p(\theta,a,b)$ is defined by~(\ref{f1_12}).
The weight $p(\theta,a,b)$ is related to the complex-valued weight $\rho(z,a,b)$ by relation~(\ref{f1_15}). 
\end{theorem}
\textbf{Proof. }
Let $a>1$ and $b\in(-x_0,x_0)\backslash\{ 0 \}$ be some fixed numbers. Consider the following function:
\begin{equation}
\label{f2_22}
g(z) = g(z,a,b) := -1/2 + {}_1 F_1 (1;a;-bz),\qquad z\in\mathbb{C}.
\end{equation}
The function $g(z)$ is an entire function and it has the following integral representation~(see~\cite[Vol. 1, p. 255]{cit_50_Bateman}):
\begin{equation}
\label{f2_24}
g(z,a,b) = -1/2 + (a-1) \int_0^1 e^{-bzu} (1-u)^{a-2} du,\qquad z\in\mathbb{C}.
\end{equation}
Choose an arbitrary $\theta\in[0,2\pi]$. We may write
$$ g(e^{i\theta}) + g(e^{-i\theta}) = -1 + 
(a-1) \int_0^1 \left(
e^{-bu e^{i\theta}} + e^{-bu e^{-i\theta}} 
\right)
(1-u)^{a-2} du = $$
\begin{equation}
\label{f2_25}
=  -1 + 2(a-1) \int_0^1 e^{-bu\cos\theta} \cos(bu\sin\theta) (1-u)^{a-2} du = 2\pi p(\theta,a,b). 
\end{equation}
This ensures that $p(\theta,a,b)$ is continuous.
Using the hypergeometric series representations we have:
$$ g(e^{i\theta}) + g(e^{-i\theta}) = -1 + 
\sum_{j=0}^\infty \frac{ (-b)^j }{ (a)_j } e^{ij\theta} + \sum_{k=0}^\infty \frac{ (-b)^k }{ (a)_k } e^{-ik\theta}. $$
Then
\begin{equation}
\label{f2_26}
p(\theta,a,b) = \frac{1}{2\pi} 
\left[
-1 + 
\sum_{j=0}^\infty \frac{ (-b)^j }{ (a)_j } e^{ij\theta} + \sum_{k=0}^\infty \frac{ (-b)^k }{ (a)_k } e^{-ik\theta}
\right].
\end{equation}
Notice that \textit{the series in~(\ref{f2_26}) converge uniformly on the interval $[0,2\pi]$}. 
This allows us to calculate the moments $m_n$ 
of $p(\theta)$:
$$ m_n := \int_0^{2\pi} e^{in\theta} p(\theta,a,b) d\theta = $$
$$ = \frac{1}{2\pi} \int_0^{2\pi} \lim_{N\rightarrow\infty} \left(
e^{in\theta}  
\left[
-1 + 
\sum_{j=0}^N \frac{ (-b)^j }{ (a)_j } e^{ij\theta} + \sum_{k=0}^N \frac{ (-b)^k }{ (a)_k } e^{-ik\theta}
\right]
\right)
d\theta = $$
$$ = \frac{1}{2\pi}
\lim_{N\rightarrow\infty} \int_0^{2\pi} e^{in\theta}  
\left[
-1 + 
\sum_{j=0}^N \frac{ (-b)^j }{ (a)_j } e^{ij\theta} + \sum_{k=0}^N \frac{ (-b)^k }{ (a)_k } e^{-ik\theta}
\right] 
d\theta = $$
\begin{equation}
\label{f2_28}
= \frac{ (-b)^n }{ (a)_n },\qquad n\in\mathbb{Z}_+.
\end{equation}
Denote
$$ S(w) := \int_\mathbb{T} w(z) \rho(z) dz,\qquad w\in\mathbb{P}, $$
where $\rho$ is from~(\ref{f1_6}), and
$$ \widetilde S(w) := \int_0^{2\pi} w(e^{i\theta}) p(\theta) d\theta,\qquad w\in\mathbb{P}. $$
The moments of the functional $S$ are known (see~\cite{cit_5100_Krall_Frink}):
$$ m_n' := S(z^n) = \frac{ \Gamma(a) }{ \Gamma(a+n) } (-b)^{n+1} = \frac{ (-b)^{n+1} }{ (a)_n }. $$
Comparing the moments $m_n$ and $m_n'$ we conclude that
\begin{equation}
\label{f2_30}
\widetilde S(w) = -\frac{1}{b} S(w),\qquad w\in\mathbb{P}. 
\end{equation}
Then
$$ \widetilde S(y_n(x,a,b) y_m(x,a,b)) = -\frac{1}{b} S(y_n(x,a,b) y_m(x,a,b)) = -\frac{c_n}{b} \delta_{n,m}, $$
and we obtain the orthogonality relations~(\ref{f1_14}).
The positivity of the weight $p$ was proved before the statement of the theorem.

Since
$$ 2\pi i \rho(z) = a-1 - \frac{b}{z} {}_1 F_1 (1;a;-b/z) = a - 1 - \frac{b}{z} \left(
g(1/z) + \frac{1}{2}
\right), $$
then
$$ g(u) = -\frac{1}{2} + \frac{a-1}{ b u } - \frac{2\pi i}{ b u } \rho (1/u),\qquad u\in\mathbb{C}\backslash\{ 0 \}. $$
Using the latter equality and relation~(\ref{f2_25}) we obtain~(\ref{f1_15}).
$\Box$

The generalized Bessel polynomials $y_n(x) = y_n(x,a,b)$ satisfy the following recurrence relation~(\cite{cit_5100_Krall_Frink}):
$$ (n+a-1)(2n+a-2) y_{n+1} = $$
\begin{equation}
\label{f2_32}
= \left[
(2n+a)(2n+a-2) \left( \frac{x}{b} \right) + a-2
\right]
(2n+a-1) y_n + n(2n+a) y_{n-1},\quad n\in\mathbb{Z}_+. 
\end{equation}

It was noticed by Beckermann in~\cite[p. 32]{cit_2100_Beckermann_2001} that one can rescale a recurrence relation
corresponding to a three-diagonal complex (not necessarily symmetric) matrix to get a recurrence relation
for a complex Jacobi matrix. We shall also apply such a normalization. 

Suppose that $a>1$, and $b\in\mathbb{R}\backslash\{ 0 \}$.
Define
\begin{equation}
\label{f2_34}
p_n(x) = p_n(x,a,b) := (-i)^n \sqrt{ \frac{ (a-1)_n (2n+a-1) }{ n! (a-1) } } y_n(x,a,b),\quad n\in\mathbb{Z}_+. 
\end{equation}

Normalized Bessel polynomials $p_n(x,a,b)$ satisfy the following recurrence relation:
\begin{equation}
\label{f2_36}
x p_n = a_{n-1} p_{n-1} + b_n p_n + a_n p_{n+1},\quad n\in\mathbb{Z}_+. 
\end{equation}
Here
\begin{equation}
\label{f2_38}
a_n := \frac{ bi }{ (2n+a) } \sqrt{ \frac{ (n+1) (n+a-1) }{ (2n+a+1) (2n+a-1) } },
\end{equation}
\begin{equation}
\label{f2_40}
b_n := -\frac{ (a-2)b }{ (2n+a) (2n+a-2) },\quad n\in\mathbb{Z}_+;\ a_{-1}:= 0,\ p_{-1}:=0,   
\end{equation}
and in the case $a=2$ one should use the limit value in~(\ref{f2_40}) for $n=0$: $b_0 = -b/2$.
Denote by $J$ the corresponding complex Jacobi matrix:
\begin{equation}
\label{f2_42}
J =
\left(
\begin{array}{ccccc}
b_0 & a_0 & 0 & 0 & \ldots \\
a_0 & b_1 & a_1 & 0 & \ldots \\
0 & a_1 & b_2 & a_2 & \ldots \\
\vdots & \vdots & \vdots & \vdots & \ddots 
\end{array}
\right),
\end{equation}
and
\begin{equation}
\label{f2_44}
\vec p(x) := (p_0(x),p_1(x),p_2(x),...)^t,
\end{equation}
where the superscript $t$ means the transposition.
Then we have
\begin{equation}
\label{f2_46}
J \vec p(x) = x \vec p(x).
\end{equation}
Denote by $\mathcal{S}$ the corresponding linear functional satisfying
\begin{equation}
\label{f2_47}
\mathcal{S} (p_n p_m) =\delta_{n,m},\qquad n,m\in\mathbb{Z}_+.
\end{equation}

Since $|a_n|$ and $|b_n|$ are uniformly bounded, the matrix $J$ defines a unique bounded operator $\mathcal{J}$ on the space $l^2$ of
square-summable complex sequences~(\cite{cit_2100_Beckermann_2001}).
We may also write~(\cite{cit_2100_Beckermann_2001}):
\begin{equation}
\label{f2_48}
\| \mathcal{J} \| \leq \mathrm{sup }_{n\geq 0} ( |a_{n-1}| + |b_n| + |a_n| ),\qquad a_{-1}=0. 
\end{equation}
Since
\begin{equation}
\label{f2_50}
|a_n| < |b| M_1,\quad |b_n| < |b| M_2,\qquad \mbox{for some } M_1,M_2>0,
\end{equation}
then
$$ \| \mathcal{J} \| \leq \mathrm{sup }_{n\geq 0} ( |a_{n-1}| + |b_n| + |a_n| ) \leq |b|(2M_1+M_2). $$
Therefore
\begin{equation}
\label{f2_52}
\| \mathcal{J} \| \leq 1,
\end{equation}
for sufficiently small $b$'s.
If~(\ref{f2_52}) holds then one can apply Theorem~3 in~\cite{cit_3000_Zagorodnyuk_2023_complex_Jacobi_matrices}
to obtain the following representation:
\begin{equation}
\label{f2_54}
\mathcal{S}(w) = \int_\mathbb{T} w(z) d\mu(z),\qquad  w\in\mathbb{P},
\end{equation}
where $\mu=\mu_{a,b}$ is a positive Borel measure on $\mathbb{T}$.
Thus, in this case we have
\begin{equation}
\label{f2_55}
\int_\mathbb{T} p_n(z,a,b) p_m(z,a,b) d\mu_{a,b}(z) =\delta_{n,m},\qquad n,m\in\mathbb{Z}_+.
\end{equation}
This fact was our motivation for studying the corresponding orthogonality relations. We could conjecture the existence
of a positive weight and study its properties.

We now turn to some applications of the orthogonality relations~(\ref{f1_14}).

\noindent
\textit{ 1) Estimates for a general solution of difference equation~(\ref{f2_36}) with $a_n,b_n$ as in~(\ref{f2_38}),(\ref{f2_40}).
}
At first we shall obtain a direct estimate for $p_n$. 

\begin{lemma}
\label{l2_1}
Let $a>1$ and $b\in\mathbb{C}\backslash\{ 0 \}$ be some fixed numbers. For an arbitrary $z\in\mathbb{C}$: $|z|>|b|$,
the following inequalities hold:
\begin{equation}
\label{f2_57}
| y_n(z,a,b) | \leq |z|^n \frac{ \Gamma(2n+a-1) }{ |b|^n \Gamma(n+a-1) } e,\qquad n\in\mathbb{N};
\end{equation}
and
\begin{equation}
\label{f2_59}
| p_n(z,a,b) | \leq |z|^n \frac{ \Gamma(2n+a-1) }{ |b|^n \Gamma(n+a-1) } 
\sqrt{ \frac{ (a-1)_n (2n+a-1) }{ n! (a-1) } }
e,\qquad n\in\mathbb{N}.
\end{equation}
\end{lemma}
\textbf{Proof. }
In conditions of the lemma we may write:
$$ |y_n(z,a,b)| = \left| 
{}_2 F_0 \left(-n,n+a-1;-;-\frac{x}{b} \right)
\right|
= \left| 
\sum_{j=0}^n \frac{ (-n)_j (n+a-1)_j (-1)^j }{ n! b^j } z^j 
\right| = $$
$$ = \left| 
\sum_{j=0}^n \frac{ (n+a-1)_j }{ (n-j)! b^j } z^j 
\right| \leq \sum_{j=0}^n \frac{ (n+a-1)_j }{ (n-j)! |b|^j } |z|^j =
\frac{ |z|^n }{ |b|^n } \sum_{j=0}^n \frac{ (n+a-1)_j }{ (n-j)! } \frac{ |b|^{n-j} }{ |z|^{n-j} } \leq
$$
$$ \leq
\frac{ |z|^n }{ |b|^n } \sum_{j=0}^n \frac{ (n+a-1)_j }{ (n-j)! } = 
\frac{ |z|^n }{ |b|^n } \sum_{k=0}^n \frac{ (n+a-1)_{n-k} }{ k! } \leq 
$$
$$ \leq
\frac{ |z|^n }{ |b|^n } (n+a-1)(n+a)...(2n+a-2)  \sum_{k=0}^n \frac{ 1 }{ k! } \leq 
\frac{ |z|^n }{ |b|^n } \frac{ \Gamma(2n+a-1) }{ \Gamma(n+a-1) }  \sum_{k=0}^\infty \frac{ 1 }{ k! }, $$
and we obtain relation~(\ref{f2_57}).
Relation~(\ref{f2_57}) follows by the definition of $p_n$.
$\Box$

Suppose that $a>1$, and $b\in\mathbb{R}\backslash\{ 0 \}$.
Denote
\begin{equation}
\label{f2_60}
q_n(z) = q_n(z,a,b) := \mathcal{S}_x \left( \frac{ p_n(z,a,b) - p_n(x,a,b) }{ z - x } \right),\quad n\in\mathbb{Z}_+, 
\end{equation}
where the subscript $x$ means that $\mathcal{S}$ acts on a polynomial in $x$.
By Theorem~\ref{t2_1} if $b\in[-x_0,x_0]$  then we have the following representation:
\begin{equation}
\label{f2_61}
q_n(z,a,b) = \int_0^{2\pi} \frac{ p_n(z,a,b) - p_n(e^{i\theta},a,b) }{ z - e^{i\theta} }
p(\theta,a,b) d\theta,\quad n\in\mathbb{Z}_+. 
\end{equation}
It is directly checked that $q_n(z,a,b)$ is a solution of the following difference equation:
\begin{equation}
\label{f2_63}
x y_n = a_{n-1} y_{n-1} + b_n y_n + a_n y_{n+1},\quad n\in\mathbb{N}. 
\end{equation}
Since $q_0 = 0$, $q_1 = 1/a_0\not=0$,
then $q_n$ and $p_n$ form a basic set of solutions of the latter difference equation.
These two sequences play an important role in the spectral theory of the corresponding Jacobi operator $\mathcal{J}$,
both in the real and the complex cases~\cite{cit_3000_Chihara}.
As in the case of real Jacobi matrices, polynomials $q_n$ can be called the polynomials of the second kind.

Let $a>1$ and $b\in\mathbb{D}\backslash\{ 0 \}$ be some fixed numbers. 
Denote
$$ k_n(z,\theta) = k_n(z,\theta,a,b) := \frac{ p_n(z,a,b) - p_n(e^{i\theta},a,b) }{ z - e^{i\theta} },\qquad  $$
\begin{equation}
\label{f2_65}
z\in\mathbb{C}\backslash\mathbb{T},\ 
\theta\in[0,2\pi],\ n\in\mathbb{Z}_+. 
\end{equation}

For arbitrary $\theta\in[0,2\pi]$ and $z\in\mathbb{C}$: $|z|>|b|$, $|z|\not=1$,
the following inequalities hold:
$$ | k_n(z,\theta,a,b) | \leq \frac{1}{ |1-|z|| } \left( |p_n(z,a,b)| + \left| p_n(e^{i\theta},a,b) \right| \right) \leq $$
\begin{equation}
\label{f2_67}
\leq \frac{ (|z|^n + 1) }{ |1-|z|| }   
\frac{ \Gamma(2n+a-1) }{ |b|^n \Gamma(n+a-1) } 
\sqrt{ \frac{ (a-1)_n (2n+a-1) }{ n! (a-1) } } e,\quad n\in\mathbb{N}.
\end{equation}

\begin{proposition}
\label{p2_1}
Let $a>1$ and $b\in(-1/3,1/3)\backslash\{ 0 \}$ be some fixed numbers. 
For an arbitrary $z\in\mathbb{C}$: $|z|>|b|$, $|z|\not=1$,
the following inequality holds:
\begin{equation}
\label{f2_69}
|q_n(z,a,b)| \leq 
\frac{ (|z|^n + 1) }{ |1-|z|| }   
\frac{ \Gamma(2n+a-1) }{ |b|^n \Gamma(n+a-1) } 
\sqrt{ \frac{ (a-1)_n (2n+a-1) }{ n! (a-1) } } e,\quad n\in\mathbb{N}.
\end{equation}
Let $\alpha,\beta\in\mathbb{C}$ be arbitrary. The general solution
\begin{equation}
\label{f2_71}
u_n(z,a,b) := \alpha p_n(z,a,b) + \beta q_n(z,a,b),\quad n\in\mathbb{N},
\end{equation}
of difference equation~(\ref{f2_63}) satisfies the following inequality:
$$ |u_n(z,a,b)| \leq $$
$$ \leq \left( |\alpha| |z|^n + |\beta| \frac{ (|z|^n + 1) }{ |1-|z|| } \right)  
\frac{ \Gamma(2n+a-1) }{ |b|^n \Gamma(n+a-1) } 
\sqrt{ \frac{ (a-1)_n (2n+a-1) }{ n! (a-1) } } e, $$
\begin{equation}
\label{f2_73}
n\in\mathbb{N}.
\end{equation}
\end{proposition}
\textbf{Proof. } 
In conditions of the proposition by~(\ref{f2_61}) we obtain that
$$ |q_n(z,a,b)| \leq \int_0^{2\pi} \left| \frac{ p_n(z,a,b) - p_n(e^{i\theta},a,b) }{ z - e^{i\theta} } \right|
p(\theta,a,b) d\theta \leq $$
\begin{equation}
\label{f2_72}
\leq \mathrm{max }_{\theta\in[0,2\pi]} | k_n(z,\theta,a,b) |.
\end{equation}
Now we can use~(\ref{f2_67}) to get~(\ref{f2_69}).
Relation~(\ref{f2_73}) follows from~Lemma~\ref{l2_1} and~(\ref{f2_69}).
$\Box$

On the other hand, we can use another representation for $q_n$. 
Let $a>1$, and $b\in\mathbb{R}\backslash\{ 0 \}$.
It is clear from the definition of $\mathcal{S}$ and orthogonality relations~(\ref{f1_5}) that
\begin{equation}
\label{f2_75}
\mathcal{S}(w) = -\frac{1}{b} \int_\mathbb{T} w(z) \rho(z,a,b) dz,\qquad w\in\mathbb{P}.
\end{equation}
In fact, the values of both sides of equality~(\ref{f2_75}) coincide on $y_n(z,a,b)$, $n\in\mathbb{Z}_+$.
Then for $z\in\mathbb{C}\backslash\mathbb{T}$ we have
\begin{equation}
\label{f2_77}
q_n(z,a,b) = -\frac{1}{b} \int_\mathbb{T} \frac{ p_n(z,a,b) - p_n(x,a,b) }{ z - x } \rho(x,a,b) dx,\quad n\in\mathbb{Z}_+. 
\end{equation}
We may write
$$ |q_n(z,a,b)| \leq \frac{1}{|b|} \int_0^{2\pi} \left| \frac{ p_n(z,a,b) - p_n(e^{i\theta},a,b) }{ z - e^{i\theta} } \right|
\left| \rho\left( e^{i\theta},a,b \right) \right|  d\theta\leq
$$
\begin{equation}
\label{f2_79}
\leq M_{a,b} \mathrm{max }_{\theta\in[0,2\pi]} | k_n(z,\theta,a,b) |, 
\end{equation}
where
\begin{equation}
\label{f2_81}
M_{a,b} := \frac{1}{|b|} \int_0^{2\pi} \left| \rho\left( e^{i\theta},a,b \right) \right| d\theta.
\end{equation}
Suppose now that $a=2$. Then
$$ \rho(z,2,b) = \frac{1}{2\pi i} e^{-b/z},\qquad z\in\mathbb{C}\backslash\{ 0 \}, $$
and
\begin{equation}
\label{f2_82}
M_{2,b} = \frac{1}{ 2\pi |b| } \int_0^{2\pi} e^{ -b\cos\theta } d\theta.
\end{equation}
Suppose additionally that $b\in[-1,1]$. Then
$$ e^{-b\cos\theta} \geq e^{-1},\qquad \theta\in[0,2\pi], $$
and
\begin{equation}
\label{f2_82_1}
M_{2,b} \geq \frac{1}{ 2\pi |b| } \int_0^{2\pi} e^{ -1 } d\theta = \frac{1}{ e |b| }.
\end{equation}
Thus $M_{2,b}$ takes large values for small values of $b$.
Comparing relations~(\ref{f2_72}) and~(\ref{f2_79}) we conclude that estimate~(\ref{f2_72}) is
better than estimate~(\ref{f2_79}) for small values of $b$.

\noindent
\textit{ 2) Similarity questions for the operator $\mathcal{J}$.
}
Denote by $\mathcal{J}_0$ the restriction of the operator $\mathcal{J}$ to the linear manifold $l^2_{fin}$.
The operator $\mathcal{J}_0$ is than similar in a generalized sence to a unitary operator in
the space~$L^2_p$. For details we refer to~\cite{cit_3000_Zagorodnyuk_2022_Similarity}.

}

\noindent
Address:

V. N. Karazin Kharkiv National University \newline\indent
School of Mathematics and Computer Sciences \newline\indent
Department of Higher Mathematics and Informatics \newline\indent
Svobody Square 4, 61022, Kharkiv, Ukraine

Sergey.M.Zagorodnyuk@gmail.com; zagorodnyuk@karazin.ua

\end{document}